\documentclass[12pt]{amsart}
\usepackage{amsfonts}
\usepackage{amssymb}
\theoremstyle{definition}

\theoremstyle{remark}

\begin{document}

\centerline{\bf ANNUAIRE \ DE \ L'UNIVERSIT\'E \ DE \ SOFIA \ "KLIMENT \ OHRIDSKI"}

\vspace{0.2in}
\centerline{\bf FACULT\'E DE MATH\'EMATIQUES ET M\'ECANIQUE}

\vspace{0.2in}
\centerline{\it Tome 75  \qquad\qquad\qquad\qquad\qquad Livre 1 -- Math\'ematiques   \qquad\qquad\qquad\qquad\qquad   1981}

\vspace{0.6in}

\vspace{0.2in}

\centerline{\large\bf THE AXIOM OF $\theta$-HOLOMORPHIC 2-PLANES }
\centerline{\large\bf IN THE ALMOST HERMITIAN GEOMETRY}

\vspace{0.3in}
\centerline{\bf Grozjo Stanilov, Ognian Kassabov}

\vspace{0.3in}
{\sl \qquad \ \ The axiom of $\theta$-holomorphic 2-planes is introduced. It is proved, 
that if an almost 

\noindent\qquad Hermitian manifold satisfies this axiom for a fixed
$\theta \in (0,\pi/2)$, then it is a real space 

\noindent\qquad form.}

\vspace{0.2in}

Let $N$ be an $n$-dimensional submanifold of a $2m$-dimensional almost Hermitian manifold $M$
with Riemannian metric $g$, almost complex structure $J$ and curvature tensor $R$. Let
$\widetilde\nabla$ and $\nabla$ denote the covariant differentiations on $M$ and $N$,
respectively. The second fundamental form $\sigma$ is a normal-bundle-valued
symmetric 2-form, defined by
$\sigma(X,Y)=\widetilde\nabla_XY - \nabla_XY$, where $X,\ Y \in {\mathfrak X}(N)$.
The submanifold $N$ is said to be totally
umbilical, if $\sigma(X,Y)=g(X,Y)H$, $H$ being the mean curvature vector of $N$,
i.e.  $H=(1/n){\rm trace}\, \sigma$. In particular, if $\sigma = 0$ \ N is called
a totally geodesic submanifold of $M$.  For $X\in \mathfrak X(N)$, 
$\xi \in \mathfrak X(N)^{\perp}$, we write 
$ 
	\widetilde\nabla_X\xi =- A_{\xi}X+D_X{\xi} \ ,
$
where $-A_{\xi}X$ (resp. $D_X\xi$)
denotes the tangential (resp. the normal) component of $\widetilde\nabla_X\xi$.
A normal vector field $\xi$ is said to be parallel, if $D_X\xi=0$ for each 
$X\in \mathfrak X(N)$.

An almost Hermitian manifold $M$ is said to be a K\"ahler manifold if 
$\widetilde\nabla J=0$. A Riemannian manifold (resp. a K\"ahler manifold) is
called a real space-form (resp. a complex space-form) if it is of constant 
sectional curvature (resp. of constant holomorphic sectional curvature).

An $n$-plane $\alpha$ in a tangent space $T_p(M)$, $p\in M$, i.e. an $n$-dimensional linear 
subspace $\alpha$ of  $T_p(M)$  is said to be 
holomorphic (resp. antiholomorphic) if  \ $\alpha =J \alpha$ \ 
(resp. \ $\alpha \perp J\alpha)$. 
An almost Hermitian manifold $M$ is said to satisfy the axiom of holomorphic $2n$-planes 
(resp. $2n$-spheres) if for each point $p \in M$ and for any holomorphic $2n$-plane
\ $\alpha$ \ in \ $T_p(M)$ \ there exists a totally geodesic
submanifold \ $N$ of $M$ \ (resp. a totally umbilical submanifold \ $N$ of $M$ \ with nonzero
parallel mean curvature vector) such that \ $p \in N$ and \ $T_pN=\alpha$. 
By changing the holomorphic $2n$-planes with antiholomorphic $n$-planes, we obtain
the axiom of antiholomorphic $n$-planes (resp. $n$-spheres).

The second author has proved in [2]:

\vspace{0.05in}
T\,h\,e\,o\,r\,e\,m A.  Let \ $M$ \ be a \ $2m$-dimensional almost Hermitian manifold,
$m \ge 2$. If \ $M$ \ satisfies the axiom of holomorphic \ $2n$-planes or \ 
$2n$-spheres for some fixed integer \ $n$, $1\le n <m$, it is an \ 
$RK$-manifold (i.e. $R(X,Y,Z,U)=R(JX,JY,JZ,JU)$ for all $X,\,Y,\,Z,\,U \in \mathfrak X(M)$) 
of pointwise constant holomorphic sectional curvature.

\vspace{0.05in}
T\,h\,e\,o\,r\,e\,m B.  Let  $M$  be a  $2m$-dimensional almost Hermitian manifold,
$m > 2$. If  $M$  satisfies the axiom of antiholomorphic \ $n$-planes or 
$n$-spheres for some  fixed integer $n$, $1< n \le m$, it  is a real space form, or
a complex space form.

\vspace{0.05in}
For the case of a K\"ahler manifold see e.g. [1, 3, 4].

For a 2-plane $\alpha$ in $T_p(M)$ the angle $ \sphericalangle(\alpha,J\alpha) \in [0,\pi/2]$
between $\alpha$ and $J\alpha$ is defined by
$$
	\cos\sphericalangle(\alpha,J\alpha) = |g(x,Jy)| \ ,
$$ 
where $ \{ x,y \}$ is an orthonormal basis of $\alpha$. Then $\alpha$ is
holomorphic (resp. antiholomorphic) if and only if $\sphericalangle(\alpha,J\alpha)=0$
(resp. $\sphericalangle(\alpha,J\alpha)=\pi/2$). In general, if 
$\sphericalangle(\alpha,J\alpha)=\theta$, $\alpha$ is called a
$\theta$-holomorphic 2-plane. Now we propose the next axiom: 

\vspace{0.05in}
{\it Axiom of $\theta$-holomorphic 2-plane (resp. 2-spheres).} For each point
$p \in M$ and for any $\theta$-holomorphic 2-plane $\alpha$ in $T_p(M)$ there exists 
a totally geodesic submanifold $N$ of $M$ (resp. a totally umbilical submanifold
$N$ of $M$ with nonzero parallel mean curvature vector) such that $p \in N$
and $T_p(N)=\alpha$.

\vspace{0.05in}
T\,h\,e\,o\,r\,e\,m. Let $M$ be a $2m$-dimensional almost Hermitian manifold, 
$m\ge 2$. If $M$ satisfies the axiom of $\theta$-holomorphic 2-planes or the
axiom of $\theta$-holomorphic 2-spheres for a fixed $\theta \in (0,\pi/2)$,
then $M$ is a real space-form. 

{\it Proof.} Let $p\in M$ and $x,\,y$ be arbitrary unit vectors in $T_p(M)$,
such that $x \perp y,\,Jy$. Then the 2-plane $\alpha$ with a basis
$ \{ x, Jx \cos\theta +y\sin\theta \} $ is $\theta$-holomorphic. Let $N$
be a totally umbilical submanifold of $M$ with parallel mean curvature vector, 
such that $p \in N$ and $T_p(N)=\alpha$. From the Codazzi's equatian
$$
	\{ R(X,Y)Z \}^{\perp} = (\overline\nabla_X\sigma)(Y,Z)-(\overline\nabla_Y\sigma)(X,Z)
$$
for $X,\,Y,\,Z \in \mathfrak X(N)$, where $ \{ R(X,Y)Z \}^{\perp} $ denotes the
normal component of $R(X,Y)Z$ and
$$
	(\overline\nabla_X\sigma)(Y,Z)=D_X\sigma(Y,Z)-\sigma(\nabla_XY,Z)-\sigma(Y,\nabla_XZ)
$$
it is easy to find
$$
	R(Jx \cos\theta +y\sin\theta,x,x,Jy)=0 \ ,   \leqno (1)
$$
$$
	R(Jx \cos\theta +y\sin\theta,x,x,Jx \sin\theta -y\cos\theta)=0 \ .   \leqno (2)
$$
We change $x$ by $-x$ in (1) and combining the result with (1) we derive
$$
	R(Jx,x,x,Jy)=0 \ .              \leqno (3)
$$
On the other hand, from (2) and (3) we obtain
$$
	H(x)=K(x,y) \ ,       \leqno (4)
$$
where
$$
	H(x)=R(x,Jx,Jx,x)  \, \qquad K(x,y)=R(x,y,y,x) \ .
$$
It follows from (4) that
$$
	H(x)=H(y)        \ .        \leqno (5)
$$
Let $m > 2$ and $u,\,v$ be arbitrary unit vectors in $T_p(M)$. We choose a unit
vector $x$ in $T_p(M)$, such that $x \perp u,\,Ju,\,v,\,Jv$. According to (5)
$$
	H(u)=H(x)=H(v) 
$$ 
i.e. $M$ is of pointwise constant holomorphic sectional curvature. Let $c=H(x)$.
Using (4) we conclude that $M$ is of pointwise constant antiholomorphic sectional
curvature $c$.

Now, let $\beta$ be an arbitrary 2-plane in $T_p(M)$ and let 
$\sphericalangle(\beta,J\beta)=\varphi$. Then it is easy to prove, that $\beta$ 
has an orthonormal basis $ \{ x, Jx\cos\varphi +y \sin\varphi \} $, where
$x,\,y$ are unit vectors in $T_p(M)$, $x \perp y,\,Jy$. Then the sectional
curvature of $\beta$ is
$$
	K(\beta)=R(x,Jx\cos\varphi +y \sin\varphi,Jx\cos\varphi +y \sin\varphi,x)
$$
and using (3), (4) and (5) we find $K(\beta)=c$. Now the assertion follows from
the Schur's theorem.

If $m=2$, let $T=R-c\pi_1$, where $c=H(x)$ and
$$
	\pi_1(x,y,z,u) =g(x,u)g(y,z)-g(x,z)g(y,u) \ .
$$ 
Then from (3), (4) and (5) we obtain easily $T=0$ and consequently $M$ is a
real space-form.

\vspace{0.05in}
C\,o\,r\,o\,l\,l\,a\,r\,y. Let $M$ be a $2m$-dimensional K\"ahler manifold, $m\ge 2$.
If $M$ satisfies the axiom of $\theta$-holomorphic 2-planes or the axiom of 
$\theta$-holomorphic 2-spheres for a fixed $\theta \in (0,\pi/2)$, then $M$ is flat.

\vspace{0.4in}

Received 13.03.1982

\end{document}